\documentclass[10pt,english, reqno]{amsart}

\usepackage[T1]{fontenc}
\usepackage{tikz}
\usepackage{amssymb,url,xspace,amsthm,adjustbox}

\usepackage[verbose,colorlinks,pagebackref]{hyperref}
\usepackage[alphabetic,nobysame]{amsrefs}

\usepackage{mathtools}
\usepackage{graphicx}

\newcommand{\bspm}{\left(\begin{smallmatrix}}
\newcommand{\espm}{\end{smallmatrix}\right)}

\usepackage{datetime}
\usepackage[T1]{fontenc}
\usepackage{chngcntr}

\theoremstyle{plain}
      \newtheorem{theorem}{Theorem}[section]
      \newtheorem{proposition}[theorem]{Proposition}
      
      \newtheorem{lemma}[theorem]{Lemma}
      \newtheorem{corollary}[theorem]{Corollary}
      \newtheorem{conjecture}[theorem]{Conjecture}
      \newtheorem{hypothesis}[theorem]{Hypothesis}

      \theoremstyle{definition}
      \newtheorem{definition}{Definition}

      \newtheorem{remark}[theorem]{Remark}
      \newtheorem{example}[theorem]{Example}

      \newtheorem*{conjecture*}{Conjecture}

\counterwithin{table}{subsection}

\newcommand{\norm}[1]{\left\lVert#1\right\rVert}

\newcommand{\mc}{\mathcal}
\newcommand{\mb}{\mathbb}
\newcommand{\mf}{\mathfrak}

\newcommand{\lb}{\left(}
\newcommand{\rb}{\right)}

\newcommand{\bc}{\mathbb{C}}

\newcommand{\bq}{\mathbb{Q}}
\newcommand{\br}{\mathbb{R}}
\newcommand{\bz}{\mathbb{Z}}

\newcommand{\tx}{\text}
\newcommand{\ten}{\otimes}

\usepackage{mathtools}
\usepackage{graphicx}

\usepackage{tikz}
\usetikzlibrary{matrix}
\usepackage{stmaryrd}
\usepackage{scrextend}
\usepackage{colortbl}
\usepackage{tabularx}
\usepackage{nicematrix}
\usetikzlibrary{arrows.meta}

\DeclareMathOperator{\GL}{GL}

\DeclareMathOperator{\SL}{SL}

\DeclareMathOperator{\Sp}{Sp}

\DeclareMathOperator{\SO}{SO}

\newcommand{\Lgroup}[1]{{\hskip-2 pt \,^L\hskip-1pt{#1}}}
\newcommand{\dualgroup}[1]{{\widehat{#1}}}

\DeclareMathOperator{\Hom}{Hom}

\DeclareMathOperator{\Ad}{Ad}

\newcommand{\ceq}{{\, :=\, }}
\newcommand{\tq}{{\ \vert\ }}

\newcommand{\Perv}{\mathsf{Per}}

\newcommand{\Loc}{\mathsf{Loc}}

\newcommand{\Ev}{\operatorname{\mathsf{E}\hskip-1pt\mathsf{v}}}
\newcommand{\Evs}{\operatorname{\mathsf{E}\hskip-1pt\mathsf{v}\hskip-1pt\mathsf{s}}}

\newcommand{\NEvs}{\operatorname{\mathsf{N}\hskip-1pt\mathsf{E}\hskip-1pt\mathsf{v}\hskip-1pt\mathsf{s}}}

\makeatletter
\newcommand{\labitem}[2]{
\def\@itemlabel{\textbf{#1}}
\item
\def\@currentlabel{#1}\label{#2}}
\makeatother

\newcommand{\1}{{\mathbbm{1}}}

\newcommand{\C}{\mathbb{C}}

\newcommand{\CL}{{\mathcal {L}}}

\newcommand{\ABV}{{\mbox{\raisebox{1pt}{\scalebox{0.5}{$\mathrm{ABV}$}}}}}

\newcommand{\ABVpure}{{\mbox{\raisebox{1pt}{\scalebox{0.5}{$\mathrm{ABV, pure}$}}}}}

\newcommand{\pure}{{\mbox{\raisebox{1pt}{\scalebox{0.5}{$\mathrm{pure}$}}}}}

\newcommand{\wh}{\widehat}

\newcommand{\bpm}{\begin{pmatrix}}
\newcommand{\epm}{\end{pmatrix}}

\usepackage{float}

\author[K. Balodis]{Kristaps Balodis}
\address{Department of Mathematics and Statistics, University of Calgary, 
2500 University Drive NW, 
Calgary, Alberta, 
T2N 1N4, 
Canada}
\email{kristaps.balodis1@ucalgary.ca}

\author[C. Cunningham]{Clifton Cunningham}
\address{Department of Mathematics and Statistics, University of Calgary, 
2500 University Drive NW, 
Calgary, Alberta, 
T2N 1N4, 
Canada}
\email{ccunning@ucalgary.ca}
\thanks{Clifton Cunningham's research is supported by NSERC Discovery Grant RGPIN-2020-05220. 
}

\author[A. Fiori]{Andrew Fiori}
\address{Department of Mathematics and Statistics, University of Lethbridge,
4401 University Drive,
Lethbridge, Alberta,
T1K 3M4,
Canada}
\email{andrew.fiori@uleth.ca}
\thanks{Andrew Fiori thanks and acknowledges the University of Lethbridge for their financial support as well as the support of NSERC Discovery Grant RGPIN-2020-05316.}

\author[Q. Zhang]{Qing Zhang}
\address{School of Mathematics and Statistics, Huazhong University of Science and Technology, Wuhan, 430074, China}
\email{qingzh@hust.edu.cn}
\thanks{Qing Zhang's research is supported by NSFC grant 12371010.}

\usepackage{datetime}
\usepackage[T1]{fontenc}
\usepackage{chngcntr}

\usepackage{amssymb}
\usepackage{mathrsfs,stmaryrd}
\usepackage{yfonts, bbm}
\usepackage{enumitem}

\usepackage{hyperref}
\hypersetup{
  colorlinks   = true, 
  urlcolor     = blue, 
  linkcolor    = blue, 
  citecolor   = green 
}

\usepackage{tikz}
\usetikzlibrary{shapes,arrows,calc,matrix}
\usepackage{tikz-cd}
\usepackage{todonotes}
\usepackage{color}

\usepackage{amsmath}
\usepackage{lipsum}
\usepackage{setspace}

\counterwithin{table}{subsection}

\makeindex

\title
[Representations and orbits with smooth closure]{Generic representations and orbits with smooth closure}

\date{\today}                                           

\begin{document}

\begin{abstract}
Let $G$ be a reductive $p$-adic group for which the local Langlands correspondence is known, and $\lambda$ an infinitesimal parameter of $G$. 
In this paper, we prove that if the $p$-adic Kazhdan-Lusztig hypothesis holds for $\lambda$, then for a Langlands parameter $\phi$ with infinitesimal parameter $\lambda$, if $\Pi_\phi(G)$ contains a generic representation, then $L(s, \phi, \Ad)$ is regular at $s=1$. 
We then prove an analogous statement for ABV-packets, which together with Vogan's conjecture on ABV-packets implies that if Arthur's conjectures for $G$ are known, then one direction of Shahidi's enhanced genericity conjecture holds:  If an Arthur packet $\Pi_\psi(G)$ contains a generic representation, then $\phi_\psi$ is tempered.  
We also offer some speculation about the relationship between Arthur type representations and singularities in varieties of Langlands parameters defined by Vogan.  
\end{abstract}

\maketitle

\setcounter{section}{-1}

\section{Introduction}

Let $G$ be a reductive algebraic group defined over a finite extension $F/\bq_p$. 
In this paper, we prove the forward direction of the following conjecture of Gross-Prasad under the assumption of (part of the desiderata of) the local Langlands correspondence and the $p$-adic Kazhdan-Lusztig Hypothesis, as well as the converse for unramified (trivial on the inertia subgroup of the Weil group) parameters:
\begin{conjecture}[\cite{GP}*{Conjecture 2.6}]\label{conj-GP}
An L-packet $\Pi_\phi(G)$ of a reductive $p$-adic group $G$ contains a generic representation if and only if $L(s, \phi, \Ad)$ is regular at $s=1$.  
\end{conjecture}
By \cite{CDFZ1}*{Proposition 3.5}, $L(s, \phi, \Ad)$ is regular at $s=1$ if and only if $C_\phi$ is open in $V_{\lambda_\phi}$ (see Section \ref{sec:KL}), and thus Conjecture \ref{conj-GP} is equivalent to the following:
\begin{conjecture}[\cite{CDFZ1}*{Conjecture 4.1}]\label{conj: genericL intro}
An L-packet $\Pi_\phi(G)$ of a reductive $p$-adic group $G$ contains a generic representation if and only if $C_\phi\subseteq V_{\lambda_\phi}$ is open.  
\end{conjecture} 
More specifically, the main result is the following.  
\begin{theorem}[Theorem \ref{thm: GP}]
Let $\lambda$ be an infinitesimal parameter of a quasi-split reductive group $G$ for which the $p$-adic Kazhdan-Lusztig hypothesis holds.   
If $\phi$ is a Langlands parameter with infinitesimal parameter $\lambda$ such that $\Pi_\phi(G)$ contains a generic representation, then $C_\phi$ is open. 

If \ $\Pi_\lambda(G)$ contains a generic representation, then if $C_\phi$ is open, $\Pi_\phi(G)$ contains a generic representation. 
\end{theorem}

As a consequence of Theorem \ref{thm: GP}, we are also able to demonstrate an analogue of Conjecture \ref{conj: genericL intro} for ABV-packets:  

\begin{corollary}[Corollary \ref{cor: genABV}]
Let $G$ be a quasi-split reductive group over a $p$-adic field $F$, and let 
    $\lambda$ be an infinitesimal parameter for which the $p$-adic Kazhdan-Lusztig hypothesis holds. 
    If \ $\Pi_\phi^\ABVpure(G)$ contains a generic representation, then $C_\phi$ is open.  

    If \ $\Pi_\lambda(G)$ contains a generic representation and $C_\phi$  is open, then $\Pi_\phi^{\ABVpure}(G)$ contains a generic representation. 
    \end{corollary}

Suppose that, in addition to the above, $G$ is a group for which Arthur's conjectures on local A-packets, as stated in \cite{Arthur:book}, are known.  If one inclusion of Vogan's conjecture (see Conjecture \ref{con: vogans conjecture}) holds, and $\psi$ is an Arthur parameter such that the $p$-adic Kazhdan-Lusztig hypothesis is known for $\lambda_\psi$.  
Then the ABV-packet analogue implies one direction of Shahidi's enhanced genericity conjecture:
\begin{conjecture}[\cite{LLS}*{Conjecture 1.2}]\label{con: enhanced shahidi intro}
An Arthur packet $\Pi_\psi(G)$ contains a generic representation if and only if $\psi$ is tempered.  
\end{conjecture}

For symplectic and split odd orthogonal groups, a proof of Conjecture \ref{con: enhanced shahidi intro} was announced in \cite{HLL} and also proved in \cite{HLLZ}.

In \cite{GI}*{Section B.2}, Gan and Ichino proved Conjecture \ref{conj-GP} under certain additional assumptions on the local Langlands correspondence (which are known to hold for classical groups and $\GL_n$). The ABV version of Shahidi's enhanced genericity conjecture is proved for classical groups in \cite{CDFZ1} using \cite{GI}*{Section B.2}. The arguments in this paper require weaker assumptions on the local Langlands correspondence than \cite{GI}*{Section B.2}, but stronger assumptions in the sense that they depend on the $p$-adic Kazhdan-Lusztig hypothesis.  

We noticed that Lemma~\ref{lem:irred_res} and Lemma \ref{lemma:memoirs}, which we use in the proof of Theorem \ref{thm: GP}, do not actually require $C_\phi$ to be open, but only use the fact that $\bar{C}_\phi$ is smooth. 
This led us to explore the broader impact of these lemmas and to consider representation-theoretic interpretations of orbits with smooth closure.  
For example, we discovered that for any orbit $C$ for which the closure $\bar{C}$ is smooth, given any irreducible equivariant local system  
$\mc{L}$ on $\bar{C}$, the representation $\pi(C, \mc{L}|_C)$ (see Section \ref{sec:KL}) appears in a single ABV packet.
We also highlight some examples supporting the speculation that for a simple group $G$, a representation $\pi(C, \mc{L})$ will only be of Arthur type if $C$ is open, closed, or $\bar{C}$ is singular. 

We now review the layout of the paper. 
In Section~\ref{background} we review background material for the paper, including some of the basic desiderata of the LLC, generic representations, standard representations, and the $p$-adic Kazhdan-Lusztig hypothesis.

In Section \ref{sec:smoothABV}, we prove that for an orbit $C$ in a variety of Langlands parameters, for which the closure $\bar{C}$ is smooth, given any irreducible equivariant local system  
$\mc{L}$ on $\bar{C}$, the representation $\pi(C, \mc{L}|_C)$ appears in a single ABV packet.  We offer some speculation in Section \ref{sec:smooth closure} about the relationship between representations of Arthur type and orbits with smooth closure by highlighting various examples.

We assume the $p$-adic Kazhdan-Lusztig hypothesis from Section \ref{sec:pklh} onward and establish some basic results about orbits with smooth closure under this assumption in Section \ref{prop:maxform}.  
The results about orbits having smooth closure are then used in Section \ref{sec:genL} and Section \ref{sec:genABV} to prove certain special cases (namely, when the $p$-adic Kazhdan-Lusztig hypothesis is known) of Conjecture \ref{conj: genericL intro}, its ABV-packet analogue, and Shahidi's enhanced genericity conjecture hold.

\subsection{Relation to other work}

Some results in this paper were reported in   \cite{CDFZ1}.

\subsection{Acknowledgments}

Kristaps Balodis would like to thank Jose Cruz, Chi-Heng Lo, Mishty Ray, and James Steele for fruitful discussions related to this paper. Also, special thanks to Maarten Solleveld for pointing out some errors in a previous version.   


\tableofcontents

\section{Background}\label{background}

\subsection{The local Langlands correspondence}\label{llc}

In this subsection, we briefly review certain basic desiderata of the local Langlands correspondence.  
In particular, we include assumptions about pure inner forms as proposed by Vogan \cite{Vogan:Langlands}, as well as assumptions on generic representations.

First, we recall the definition of generic representations.  Let $F$ be a non-archimedean local field, $G$ a connected quasi-split reductive group, and $B$ be a Borel subgroup of $G$ with unipotent radical $U$. The torus $T=B/U$ acts on $\Hom(U,\C^\times)$. A character $\theta: U(F)\to \C^\times$ is called \emph{generic character}\index{generic character} if its stabilizer in $T(F)$ is the center $Z(F)$ of $G(F)$.  
For a generic character $\theta$ of $U(F)$, a representation $\pi$ of $G(F)$ is called $\theta$-generic\index{generic, representation} if $\Hom_{U(F)}(\pi,\theta)\ne 0$; this depends only on the $T(F)$-orbit of $\theta$. 
If $\theta$ is understood from the context, we simply say $\pi$ is generic.

 A \emph{Whittaker datum} \index{Whittaker datum} $\mathfrak{w}$ \index{$\mf{w}$}for $G$ is a $G(F)$-conjugacy class of pairs $(B,\theta)$, where $B$ is a Borel subgroup of $G$ and $\theta$ is a generic character of $U(F)$, where $U$ is the unipotent radical of $B$. 

We now suppose $G$ is a connected reductive algebraic group (not necessarily quasi-split) over a non-Archimedean field $F$.  As in \cite{Vogan:Langlands}*{Definition 2.6, Proposition 2.7} we call each element $\delta$ of the non-commutative cohomology set $H^1(F, G)$, a \emph{pure inner form}\index{pure inner form} of $G$.
Each $\delta$\index{$\delta$} determines a form $G_\delta$ of $G$. For each pure inner form $\delta$, let $\Pi(G_\delta)$ be the set of isomorphism classes of  smooth irreducible representations of $G_\delta(F)$ and set 
\[
\Pi^{\pure}(G)\ceq\coprod_{[\delta]\in H^1(F,G)}\Pi(G_\delta).
\]
\index{$\Pi^{\pure}$}
Let $\hat{G}$ be the complex dual group, let $  W_F$ be the Weil group, and define ${}^LG=\hat{G}(\bc)\rtimes W_F$. 
We write $\Phi(G)$\index{$\Phi(G)$} for the set of equivalence classes of Langlands parameters $\phi:W_F\times\SL_2(\bc)\to {}^LG$ of $G$. For every Langlands parameter $\phi$, we define $A_\phi\ceq \pi_0(Z_{\dualgroup{G}}(\phi))$.

We now describe the desiderata of the local Langlands correspondence, which are relevant for our purposes. 
\begin{enumerate}
\item
There is a finite-to-one surjective map 
\[
\Pi^{\pure}(G)\to \Phi(G).
\]
 For each Langlands parameter $\phi$, the \emph{pure} L-packet of $\phi$ \index{pure, L-packet} 
\[
\Pi_\phi^{\pure}(G)\subset \Pi^{\pure}(G)
\]
is defined to be the pre-image of $[\phi]$. 

\item\label{desi2}  
For a fixed Whittaker datum $\mathfrak{w}$, there is a family of bijections \index{$J(\mathfrak{w})$}
\[
J(\mf{w}): \Pi_\phi^{\pure}(G)\to \wh{A_\phi},
\]
for each Langlands parameter $\phi$ of $G$, such that if $\Pi^\pure_\phi(G)$ contains a $\mf{w}$-generic representation $\pi$, then it is the unique $\mf{w}$-generic representation in $\Pi_\phi^\pure(G)$ and $J(\mathfrak{w})(\pi)$ is the trivial representation of $A_\phi$.\index{$A_\phi$}
We call $J(\mf{w})$\index{$J(\mf{w})$} a \emph{Whittaker normalization}.\index{Whittaker normalization} 
\end{enumerate}

We say that $\Pi_\phi(G)$ (resp. $ \Pi_\phi^\pure(G)$) and $\phi$ are \emph{generic}\index{generic, Langlands parameter}\index{generic, L-packet} if $\Pi_\phi(G)$ contains a generic representation (See, for example, \cite{CDFZ1}*{Section 1.5}).


\subsection{Kazhdan-Lusztig hypothesis for \texorpdfstring{$p$}{}-adic groups}\label{sec:KL}

In this subsection, we recall the notion of \emph{Vogan varieties}\index{Vogan varieties} which are the main focus of this paper.  
We will also recall the statement of the $p$-adic Kazhdan-Lusztig hypothesis, \index{$p$-KLH} which is crucial for the results in the latter half of the paper. 
In short, the $p$-adic Kazhdan-Lusztig hypothesis articulates a deep relationship between equivariant perverse sheaves on moduli spaces of Langlands parameters introduced by \cite{Vogan:Langlands} and the multiplicities of irreducible representations in \emph{standard representations}, whose definition depends on the following result:

\begin{theorem}\label{thm: Konno}{\cite{Kon}*{Theorem 3.5}}
For any irreducible admissible representation $(V, \pi)$ of $G(F)$, there exists a parabolic $P=MU$ and an irreducible tempered representation $\tau \in \Pi_{\mathrm{temp}}(M(F))$ and $\mu\in\mf{a}_P^{\ast, +}$ such that $\pi$ is the unique irreducible quotient of $I_P^G(e^\mu\ten \tau)$, where $I_P^G$ is the normalized parabolic induction functor.  
Moreover, the triple $(P, \mu, \tau)$ is unique up to conjugation by the Weyl group of $G$. 
\end{theorem}

See \cite{Kon} for more details and unexplained notations regarding the above theorem. 
	
The important take-away from the above theorem is that fixing a choice of minimal parabolic $P_0$ of $G$, for any irreducible representation $\pi$ of $G$, there is a unique representation $I_P^G(e^\mu\ten \tau)$ such that $P$ is a standard parabolic and $\pi$ is the unique irreducible quotient of $I_P^G(e^\mu\ten \tau)$.  We will refer to $I_P^G(e^\mu\ten \tau)$ as the standard representation of $\pi$ and denote it by $S_\pi$.  \index{$S_\pi$}\index{standard representation}

We now turn to the description of the moduli space of Langlands parameters mentioned at the start of this section. 
Given a Langlands parameter
$$\phi:W_F\times \SL_2(\bc)\to {}^LG,$$
its infinitesimal parameter \index{infinitesimal parameter}
$$\lambda_\phi:W_F\to {}^LG,$$
 is defined by 
$$\lambda(w)=\phi\lb w, d_w\rb,$$
where
$$d_w:=\begin{pmatrix}
    \norm{w}^{1/2} & 0 \\
    0 & \norm{w}^{-1/2}
\end{pmatrix},$$
for the norm map $\norm{\cdot}:W_F\to \br^\times$. 

Given a fixed infinitesimal parameter $\lambda$, we define 
$$V_\lambda:=\{x\in \hat{\mf{g}} \ |\  \forall w\in W_F, \Ad(\lambda(w))x=\norm{w}x\}.$$
We refer to $V_\lambda$ as a \emph{Vogan variety}\index{Vogan variety}, as it was introduced in \cite{Vogan:Langlands}. 

The group
$$H_\lambda:=\{g\in \hat{G}(\bc)\ | \ \forall w\in W_F, g\lambda(w)g^{-1}=\lambda(w)\},$$
acts on $V_\lambda$ by conjugation.
There are finitely many $H_\lambda$-orbits on $V_\lambda$, and there is a unique closed orbit and a unique open orbit.  
The closure of each orbit is a union of orbits.  
For orbits $C$ and $D$, we write $D\leq C$ if $D\subseteq \bar{C}$.  

The orbits of $V_\lambda$ are in bijection with Langlands parameters having infinitesimal parameter $\lambda$.  
In particular, if $\phi(w, d_w)=\lambda(w)$, then the restriction
$$\phi|_{\SL_2(\bc)}:\SL_2(\bc)\to {}^LG,$$
determines an algebraic map 
$$\SL_2(\bc)\to \hat{G}(\bc),$$
and thus defines an $\SL_2$-triple $(h_\phi, x_\phi, y_\phi)$ such that $x\in V_\lambda$.  
Two Langlands parameters $\phi, \phi'$ with $\lambda_\phi=\lambda=\lambda_{\phi'}$ are equivalent if and only if $x_{\phi}$ and $x_{\phi'}$ are in the same $H_\lambda$-orbit.

Recall from Desiderata \ref{desi2} that fixing a choice of Whittaker normalization provides a bijection $J(\mf{w}):\Pi_\phi^\pure(G)\to\wh{A_\phi} $ for each Langlands parameter such that if $\pi$ is $\mf{w}$-generic, then $J(\mf{w})(\pi)$ is the trivial representation. 
Letting $x\in V_\lambda$ and $H_\lambda(x)$ be the stabilizer of $x$ in $H_\lambda$, the inclusion $Z_{\hat{G}}(\phi_x)\to H_\lambda(x)$ induces an isomorphism
$$A_{x}:=H_\lambda(x)/H_\lambda(x)^\circ\cong Z_{\hat{G}}(\phi)/Z_{\hat{G}}(\phi)^\circ.$$

\begin{definition}\label{piCL}
    Given some fixed choice of Whittaker normalization $J(\mf{w})$, if $\mc{L}\in \Loc_{H_\lambda}(C)$ is irreducible, let $\pi(C, \mc{L})$\index{$\pi(C, \mc{L})$} be the inverse image of $\mc{L}$ under the composition
$$\Pi_\phi(G)\xrightarrow{J(\mf{w})} \wh{A_{\phi_x}} \to \wh{A_x}\to \Loc_{H_\lambda}(C).$$
We write $S(C, \mc{L}):=S_{\pi(C, \mc{L})}$ for the corresponding standard representation. 
\end{definition}

See \cite{CDFZ1}*{Section 2.2} for more details on making a change of Whittaker normalization, and its compatibility with the local Langlands correspondence.  

\begin{hypothesis}[$p$-adic analogue of the Kazhdan-Lusztig Hypothesis]\label{hypothesis}
Let $C, D\subseteq V_\lambda$ be $H_\lambda$-orbits of $V_\lambda$, and $\mc{L}\in \Loc_{H_\lambda}(C), \mc{F}\in \Loc_{H_\lambda}(D)$ be irreducible local systems.  
We denote by $\mc{IC}(D, \mc{F})$\index{$\mc{IC}(D, \mc{F})$} the intersection cohomology complex associated to $\mc{F}$. 
Let $\pi=\pi(C, \mc{L})$, and $\sigma=\pi(D, \mc{F})$.  
  Then, the multiplicity $[S_\pi: \sigma]$ with which $\sigma$ appears in the Jordan-H\"older series of $S_\pi$ is given by
    $$[S_\pi: \sigma] =\sum_{n\in \bz}\dim\Hom_{\Loc_{H_\lambda}(C)}\left( \mc{H}^n(\mc{IC}(D,  \mc{F}))|_{C}, \mathcal{L}\right).$$
\end{hypothesis}

The Kazhdan-Lusztig hypothesis for $p$-adic groups first appeared in \cite{Zelevinskii:KL} for general linear groups. Here we review the general form of the conjecture; see also \cite{Vogan:Langlands}*{Conjecture 8.11}, \cite{CFMMX}*{\S 10.3.3}, \cite{Solleveld:pKLH}*{Theorem E, (b) and (c)}, and, for real groups, \cite{ABV}*{Corollary 1.25}.

The result \cite{CG}*{Theorem 8.23} has occasionally been referred to as a proof of the $p$-adic Kazhdan-Lusztig hypothesis. 
 However, while the work of \cite{CG} provides substantial progress towards the $p$-adic Kazhdan-Lusztig hypothesis, their result only applies to quasi-split $p$-adic groups $G$ for which the dual group $\hat{G}$ is semi-simple and simply connected, and only for those representations of $G$ with a non-zero Iwahori-fixed vector.  
In fact, the results \cite{CG} are written entirely in terms of modules for the Iwahori-Hecke algebra of $G$, and it remains to be verified that the so-called "standard modules" of \cite{CG} actually correspond to the standard representations as defined in this document.  

The more recent result \cite{Solleveld:pKLH}*{Theorem 5.4} proves an analogue of the $p$-adic Kazhdan-Lusztig hypothesis, but for modules over affine graded Hecke algebras.  
The result of \cite{Solleveld:pKLH}*{Theorem 5.4} applies to a much wider case of groups than \cite{CG}.   
 Together with \cite{Sol2}*{Lemma 6.2}, which ensures that standard \emph{modules} of graded affine correspond to standard representations as defined in this paper, the $p$-adic Kazhdan-Lusztig hypothesis as stated by \cite{Vogan:Langlands}, holds in the cases listed in \cite{Solleveld:pKLH}*{Theorem 5.4}.  
 In particular, for various groups, it holds beyond just those representations generated by their Iwahori-fixed vectors. 

 The $p$-adic Kazhdan-Lusztig hypothesis has been verified for unipotent representations of $G_2$ in \cite{CFZ:cubic} and \cite{CFZ:unipotent}.


\subsection{\texorpdfstring{$\NEvs$}{} and ABV-packets}

Since $V_\lambda$ is a vector space, and $H_\lambda$ acts linearly on $V_\lambda$, there is an induced $H_\lambda$-action on $V_\lambda^\ast$.  
Moreover, the cotangent bundle is trivial $T^\ast(V_\lambda)\cong   V_\lambda\times V_\lambda^\ast$.  
Consider the closed subvariety
\[
\Lambda \ceq \{ (x,\xi)\in V_\lambda\times V_\lambda^* \tq [x,\xi]=0 \},
\]
of $T^\ast(V_\lambda)$ where, $[~,~]$ is the Lie bracket in $\widehat{ \mathfrak{g}}$.
Note that while $\hat{G}$ may not be semi-simple, we can still identify $V_\lambda^\ast$ with a subset of $\hat{\mf{g}}$ as in \cite{CFMMX}*{Section 6.2-6.3}.  
Likewise, the subvariety
\[
\Lambda^\ast \ceq \{ (\xi,x)\in V_\lambda^* \times V_\lambda \tq [x,\xi]=0 \}
\]
of $T^\ast(V_\lambda^\ast)\cong V_\lambda\times V_\lambda^\ast$ may be identified with $\Lambda$.  
We write $p : \Lambda \to V_\lambda$ and $q : \Lambda \to V_\lambda^\ast$ for the obvious projections, and for every pair of $H_\lambda$-orbits $C\subseteq V_\lambda, D\subseteq V_\lambda^\ast$, set $\Lambda_{C}\ceq p^{-1}(C)$ and $\Lambda_{D}^*\ceq q^{-1}(D)$.

Pyasetskii duality $C \mapsto C^*$ defines a bijection between the $H_\lambda$-orbits in $V_\lambda$ and the $H_\lambda$-orbits in $V_\lambda^*$ and is uniquely characterized by the following property: under the isomorphism $\Lambda \to \Lambda^*$ defined by $(x,\xi)\to (\xi,-x)$, the closure of $\Lambda_{C}$ in $\Lambda$ is isomorphic to the closure of $\Lambda_{C^*}^*$ in $\Lambda^*$.
This duality may also be characterized by passing to the regular conormal variety as follows.
Set
\[
    \Lambda_C^\mathrm{reg}=\{(x,\xi)\in \Lambda_C \tq \xi \in q(\Lambda_C)\setminus q(\bar{\Lambda}_{C'})\text{ where }C \subseteq \bar{C'} \text{ but } C \neq C' \}, 
\]
    where $q : \Lambda \to V_\lambda^*$ is the projection and also set
$
\Lambda^\mathrm{reg} \ceq \bigcup_C \Lambda_{C}^\mathrm{reg}
$ (this is a disjoint union, in fact).
Likewise, we define $\Lambda_{D}^{*,\mathrm{reg}}$ for an $H$-orbit $D\subset V_\lambda^*$.

In  \cite{CFMMX}*{Section 7.10}, a functor
\[
\NEvs_{C}: \Perv_{H_{\lambda}}(V_\lambda)\to \Loc_{H_{\lambda}}(\Lambda_{C}^{\mathrm{gen}}),
\]
is defined (in terms of a vanishing cycles functor) for every $H_\lambda$-orbit $C$ in $V_\lambda$, where  $\Lambda^\mathrm{gen}_{C}\subseteq \Lambda$ is a connected $H_\lambda$-stable open subset defined in \cite{CFMMX}*{Section 7.9}.

Fix a Whittaker datum $ \mathfrak{w}$ for $G$. 
By Desiderata \ref{desi2} and Definition \ref{piCL}, the representations of $\Pi_\lambda^\pure(G)$ are all of the form $\pi(C, \mc{L})$ for an irreducible $H$-equivariant local system on an $H$-orbit $C$. 

\begin{definition}
For a local Langlands parameter $\phi$ with infinitesimal parameter $\lambda$, and a Whittaker datum $\mathfrak{w}$ as above, define
\begin{equation}\label{definition:ABVw}
\Pi_{\phi,\mathfrak{w}}^{\ABVpure}(G)
:=
\{ \pi \in  \Pi_\lambda^\pure(G) \tq \pi=\pi(C, \mc{L}), \NEvs_{C_\phi}(\mc{IC}(C, \mc{L}))\ne 0\}.
\end{equation}

By \cite{CDFZ1}*{Theorem 2.2}, $\Pi_{\phi, \mf{w}}^\ABVpure(G)$ is independent of the choice of Whittaker datum $\mf{w}$, and thus we write
$$\Pi_\phi^\ABVpure(G):=\Pi_{\phi, \mf{w}}^\ABVpure(G).$$
\end{definition}

Set $W_F'':=W_F'\times \SL_2(\bc)=W_F\times \SL_2(\bc)\times\SL_2(\bc)$.  An \emph{Arthur parameter} of $G$ is a homomorphism 
$$\psi:W_F''\to {}^LG,$$
such that $\phi_\psi(w, x):=\psi(w, x, d_w)$ is a \index{$\phi_\psi$} Langlands parameter of $G$, the restriction to the second $\SL_2(\bc)$-factor induces an algebraic map $\SL_2(\bc)\to \hat{G}$, and the image of $\psi|_{W_F}$ has bounded image.
We say that a Langlands parameter $\phi$ is of Arthur type if there is an Arthur parameter $\psi$ such that $\phi$ is equivalent to $\phi_\psi$. 
We will say that an $H_\lambda$-orbit $C\subseteq V_\lambda$ is of Arthur type if there is an Arthur parameter $\psi$ such that $C=C_{\phi_\psi}$.  
We let $\Psi(G)$ denote the set of all equivalence classes of Arthur parameters under $\hat{G}$-conjugation of the image.

To each Arthur parameter $\psi$ one can associate a subset
$$\Pi_\psi(G)\subseteq \Pi(G),$$
called an \emph{Arthur packet} (or A-packet) which depends only on the equivalence class of $\psi$.  
We will say that an irreducible representation $\pi$ of $G(F)$ is of \emph{Arthur type} if there is an Arthur parameter $\psi$ for which $\pi\in \Pi_\psi(G)$, and that a Langlands parameter is of Arthur type if it is equivalent to $\phi_\psi$ for some Arthur parameter $\psi$. 
The precise definition can be found in \cite{Arthur:book}.\index{$\Pi-\psi(G)$}
We also define 
$$\Pi_\psi^\pure(G):=\coprod_{[\delta]\in H^1(F, G)}\Pi_\psi(G_\delta).$$

The following conjecture was originally posed by Vogan using the language of microlocal geometry and reformulated in terms of the $\NEvs$ functor in \cite{CDFZ1}*{Conjecture 1. a), Section 8.3}.

\begin{conjecture}\label{con: vogans conjecture}
If $G$ is a $p$-adic group for which Arthur's conjectures \cite{Arthur:Conjectures} are known, then
    $$\Pi_\psi^\pure(G)= \Pi_{\phi_\psi}^\ABVpure(G)$$
for all Arthur parameters $\psi$ for $G$.
\end{conjecture}

The proof of Vogan's conjecture for $G=\GL_n$ first appeared in \cite{CR} and \cite{CR2}.  
For unipotent representations of $\GL_n$, a combinatorial proof was given in \cite{Rid}, and another independent proof in \cite{Lo}.
This conjecture is proved for quasi-split classical groups and their pure inner forms in \cite{CHLLRXS}.


\section{Orbits with smooth closure and Arthur type representations}

In this section, let $G$ be a reductive group over a $p$-adic field. We assume Arthur's conjectures \cite{Arthur:Conjectures} are known - namely, there exists a local Arthur packet theory for $G$. 

\subsection{ABV-packets for orbits with smooth closure}\label{sec:smoothABV}

In this section, we prove Proposition \ref{prop: singleton}: If $C$ is an orbit of Arthur type, and $\mc{L}$ is a certain kind of $H_\lambda$-equivariant local system (of which the trivial local system is an example), then it belongs to a single ABV packet. 
We also explore some relationships between representations of Arthur type, that is, irreducible representations $\pi$ for which $\pi\in \Pi_\psi(G)$ for an Arthur parameter $\psi$, and orbits with smooth closure in the corresponding Vogan variety.  

Following \cite{HLLA}, given an Arthur parameter $\psi$ of $G$ and an irreducible representation $\pi\in \Pi_\psi(G)$, we define
$$\Psi(\pi)=\{[\psi]\in \Psi(G) : \pi\in \Pi_\psi(G) \}.$$
Since $A$-packets usually have nontrivial intersections, $\Psi(\pi)$ is usually not a singleton. Various structures on $\Psi(\pi)$ were studied in \cites{HLLA, HLLZ} when $G$ is a classical group.  

Shahidi's enhanced genericity conjecture \cite{LLS}*{Conjecture 1.2} predicts that $\Pi_\psi(G)$ contains a generic representation iff $\phi_\psi$ is tempered, which generalizes a previous hypothesis \cite{Shahidi: plancherel}*{Conjecture 9.4}.  
Some progress towards \cite{LLS}*{Conjecture 1.2} was made in \cite{CDFZ1}*{Proposition 4.8}. 

Inspired by Shahidi's conjecture, \cite{CDFZ1}*{Conjecture 4.3} posits that $\Pi_\phi^{\ABVpure}(G)$ contains a generic representation if and only if $C_\phi$ is open. 
Progress towards this result is made in \cite{CDFZ1}*{Theorem 4.4}.  

Recall that $\Phi(G)$ denotes the set of equivalence classes of 
Langlands parameters of $G$.  
Given an irreducible representation $\pi$ of $G$, we define
$$\Phi^{\ABV}(\pi):=\{[\phi]\in \Phi(G): \pi\in \Pi_\phi^{\ABV}(G)\}.$$

Note that \cite{CDFZ1}*{Conjecture 4.3} implies that if $\pi$ is generic, then $\Phi^{\ABV}(\pi)$ is a singleton.  
If we assume Conjecture \ref{conj: genericL intro}, then conversely, if $\Phi^\ABV(\pi)$ is a singleton for every generic representation, \cite{CDFZ1}*{Conjecture 4.3} holds.  
First we note that even without the assumption on $\Phi^\ABV(\pi)$, if $C_\phi$ is open, then by Conjecture \ref{conj: genericL intro} 
$$\Pi_\phi^\pure(G)\subseteq \Pi_\phi^\ABVpure(G),$$
contains a generic representation. 

Meanwhile, if for some $\phi$, $\Pi_\phi^\ABVpure(G)$ contains a generic representation $\pi$, then
$$\pi\in \Pi_{\phi_\pi}^\ABVpure(G)\subseteq \Pi_{\phi_\pi}^\ABVpure(G),$$
and since $\Phi^\ABV(\pi)$ is a singleton, it must be that $\phi$ is equivalent to $\phi_\pi$.  
By Conjecture \ref{conj: genericL intro} $C_{\phi_\pi}=C_\phi$ is open.

Before stating the main results of this section, it will be useful to have the following lemma:

\begin{lemma}\label{lem:irred_res}
Suppose that an $H_\lambda$-orbit $C$ of a Vogan variety $V_\lambda$ has a smooth closure. 
If $\mc{L}\in \Loc_{H_\lambda}(\bar{C})$ is irreducible, then so is $\mc{L}|_C$. 
\end{lemma} 
\begin{proof}
    Since $\bar{C}$ is smooth, $\mc{L}[\dim \bar{C}]=\mc{L}[\dim C]$ is perverse.  
    Hence, there must exist an orbit $D\leq C$ and an $H$-equivariant local system $\mc{F}$ on $D$ such that $\mc{IC}(D, \mc{F})\cong \mc{L}[\dim C]$.  
    Thus 
    $$\mc{IC}(D, \mc{F})|_C\cong \mc{L}|_C[\dim C],$$
    and thus it has non-zero support on $C$, but this is only possible if $D=C$, in which case 
    $$ \mc{L}|_C[\dim C]\cong \mc{IC}(C, \mc{F})|_C\cong \mc{F}[\dim C],$$
    and therefore $\mc{L}|_C\cong \mc{F}$, which concludes the result. 
\end{proof}

\begin{lemma}\label{lemma:memoirs}
    Let $C\subseteq V_\lambda$ be an $H_\lambda$-orbit such that $\bar{C}$ is smooth, and suppose $\mc{L}\in \Loc_{H
    _\lambda}(\bar{C})$ is irreducible.    
    Then, for all $H_\lambda$-orbits $D\subseteq \bar{C}$ such that $D\neq C$,
    $$\NEvs_D(\mc{IC}(C, \mc{L}|_C))\cong 0.$$
\end{lemma}
\begin{proof}
    The functor $\NEvs$ is defined in terms of another functor $\Ev$, defined in terms of vanishing cycles at the end of \cite{CDFZ1}*{Section 7.3}.  
    It suffices to determine when $\Ev$ is 0, which is defined
    $$(\Ev_D \mc{F})_{(x, \xi)}=(\mathrm{R}\Phi_{f_\xi}\mc{F})_x=\mathrm{R}\Phi_{f_{D^\ast}}(\mc{F}\boxtimes \mathbbm{1}_{D^\ast})_{(x, \xi)}.$$

    Given that the intermediate extension functor respects composition, we conclude that for the inclusion $j:\bar{C}\hookrightarrow V$,
    $$\mc{IC}(C, \mc{L}|_C)\cong j_!\mc{L}[\dim C].$$
    Therefore
    $$\Ev_D(\mc{IC}(C, \mc{L}|_C)\cong \mathrm{R}\Phi_{f_{D^\ast}}(j_!(\mc{L})\boxtimes \mathbbm{1}_{D^\ast})|_{\Lambda_{D}^\text{reg}}.$$
 By \cite{CFMMX}*{Proposition 7.10} we know that this is only non-zero if $D\subseteq \bar{C}$. Thus, by proper base change, it suffices to compute
 \[ \mathrm{R}\Phi_{f|_{\overline{C}\times D^*}}(\mc{L} \boxtimes \mathbbm{1}_{D^\ast})|_{\Lambda_{D}^{\mathrm{reg}}}  \]
     where we now write $f=f_{D^\ast}$.

     By \cite{CFMMX}*{Lemma 7.3}, since $\mc{L}\boxtimes \1_{D^\ast}$ is a local system and $\overline{C}\times D^*$ is smooth, $\mathrm{R}\Phi_f(\mc{L}\boxtimes \mathbbm{1}_{D^\ast})_{(x, \xi)}\cong 0$ for $(x, \xi)$ in the smooth locus of $f$.  
     We claim that $\Lambda_{D}^{\mathrm{reg}}$ is contained in the smooth locus, which would give the result.
     
Thus, it suffices to determine the singular locus $f$. Now, using that $\overline{C}\times D$ is smooth, $\bar{C}\times D^* \subset V_\lambda\times V_\lambda^*$, and $f$ is the restriction of the perfect pairing of the pairing between $V_\lambda$ and $V_\lambda^*$ to $\overline{C}\times D$, the Jacobian condition for smoothness tells us $f$ is singular at $(x,\xi)$ if and only if
    \[ df|_{(x,\xi)} \in {\rm span}\{ dg|_{(x,\xi)} \;|\; g\in I(\overline{C}\times \overline{D^*} )\}, \]
    where $I$ is the ideal of the functions defining $\overline{C}\times \overline{D^*} \subset V_\lambda\times V_\lambda^*$.
    
Since $\overline{C}\times \overline{D^*} \subset V_\lambda\times V_\lambda^*$ we may interpret all of the differentials as being in $V_\lambda^*\times V_\lambda$. From the definition of $f$, we then have $df|_{(x,\xi)} = (\xi,x)$. Moreover, as $\overline{C}\times \overline{D^*}$ is a product, we have
    \[  {\rm span}\{ dg|_{(x,\xi)} \;|\; g\in I(\overline{C}\times \overline{D^*} ) \} = 
    {\rm span}\{ dg|_{x} \;|\; g\in I(\overline{C}) \}  \oplus {\rm span}\{ dg|_{\xi} \;|\; g\in I( \overline{D^*} ) \}. \]
    Thus $(x,\xi)$ is a smooth point unless 
    \[ \xi \in   {\rm span}\{ dg|_{x} \;|\; g\in I(\overline{C}) \}  \text{ and } x \in  {\rm span}\{ dg|_{\xi} \;|\; g\in I( \overline{D^*} ) \}, \]
    but this is precisely to say that
    \[  \xi \in T_{\overline{C}}^*|_x \text{ and } x \in T_{D^*}^*|_{\xi}. \]
    The first condition gives $(x,\xi) \in \overline{T_C^*}$, and hence, by the definition of $\Lambda_{D}^{\mathrm{reg}}$, we have $(x,\xi)\not\in \Lambda_{D}^{\mathrm{reg}}$.  
\end{proof}
\begin{remark}
The second condition $x \in T_{D^*}^*|_{\xi}$ obtained above does not actually rely on the smoothness of $\overline{C}$ and, indeed, gives the conclusion of \cite{CFMMX}*{Lemma 7.21}, namely, that the support of $\mathrm{R}\Phi_{f|_{\overline{C}\times D^*}}(\mc{L}_{\overline{C}} \boxtimes \mathbbm{1}_{D^\ast}) $ should be in $\Lambda_D$.
\end{remark}

Recall that by Lemma \ref{lem:irred_res}, for an irreducible $\mc{L}\in\Loc_{H_\lambda}(\bar{C})$, the restriction $\mc{L}|_C$ is irreducible and thus $\pi(C, \mc{L}|_C)$ is well-defined.  
Hence, we obtain the following corollary.

\begin{proposition}\label{prop: singleton}
 Let $G$ be a reductive group over a $p$-adic field $F$, and assume the local Langlands correspondence for $G$.
 Let $\phi$ be a Langlands parameter with infinitesimal parameter $\lambda$, and suppose that the closure of $C=C_\phi$ is smooth.  
 If $\mc{L}$ is an irreducible $H_{\lambda_\phi}$-equivariant local system on $\bar{C}_\phi$, then
$$\Phi^{\ABV}(\pi(C,\CL|_C))=\{\phi\}.$$
\end{proposition}
\begin{proof}
    Let $\mc{L}$ be an $H$-equivariant local system on $\bar{C}_\phi$.  
 By \cite{CFMMX}*{Proposition 7.10} (and using the fact that $\NEvs$ is non-zero if and only if $\Evs$ is non-zero), $\NEvs_{C_{\phi'}}(\mc{IC}(C_\phi, \mc{L}|_{C_{\phi}}))$ is only nonzero for those Langlands parameters $\phi'$ which are contained in the support $\bar{C}_\phi$ of $\mc{IC}(C_\phi, \mc{L}|_{C_\phi})$.  

    By Lemma \ref{lemma:memoirs}, for a Langlands parameter $\phi'$ such that $C_{\phi'}\subseteq \bar{C}_\phi$, and $C_{\phi'}\neq C_\phi$, we have
    $$\NEvs_{C_{\phi'}}(\mc{IC}(C_\phi, \mc{L}|_{C_\phi}))\cong 0.$$
    Finally, by \cite{CFMMX}*{Proposition 7.13 b), Proposition 7.19}
    $$\NEvs_{C_\phi}(\mc{IC}(C_\phi, \mc{L}_{C_\phi}))\neq 0,$$
    and thus $\pi(C_\phi, \mc{L}|_{C_\phi})\in \Pi_\phi^\ABV(G)$, which concludes the result.  
\end{proof}

\begin{remark}
    Note that the assumption that the local system in $\pi(C, \mc{L}|_C)$ is the restriction of an irreducible $H_\lambda$-equivariant local system on the closure is crucial.  
    In Example \ref{ex:SO(7)}, we discuss the existence of an orbit with smooth closure $\bar{C}_3$ for which there is a non-trivial local system whose representation $\pi(\phi_3, -)$ belongs to more than one ABV-packet. 
    Moreover, the converse statement fails.  
    In \cite{CFMMX}*{Chapter 16}, one can see that there are two representations  $\pi_4^+$ and $\pi_5^+$ which belong to a single ABV-packet, and yet their orbit closures are not smooth.  
    However, we note that neither $\pi_4^+$ nor $\pi_5^+$ is of Arthur type. 
\end{remark}

A special case of Proposition \ref{prop: singleton} states that if $\bar{C}$ is smooth, then $\Phi^\ABV(\pi(C, \1))=\{\phi\}$.  

\begin{corollary}\label{cor: closed single}
     Let $\phi$ be a Langlands parameter such that $\bar{C}_\phi$ is smooth, and let $\mc{L}$ be an irreducible $H_\lambda$-equivariant local system on $\bar{C}_\phi$.  
     Assuming Vogan's Conjecture~\ref{con: vogans conjecture}, $\pi(C_\phi, \mc{L}|_{C_\phi})$ is Arthur type if and only if the Langlands parameter $\phi$ is Arthur type. 
     Moreover, in this case $\pi(C_\phi, \mc{L}|_{C_\phi})$ is contained in a unique Arthur packet.  
\end{corollary}
    \begin{proof} 
If $\bar{C}_\phi$ is smooth, and $\mc{L}$ is an irreducible $H_\lambda$-equivariant local system on $\bar{C}_\phi$, then by Proposition \ref{prop: singleton} $\Phi^{\ABV}(\pi(C_\phi, \mc{L}|_{C_\phi}))=\{\phi\}$. 

If $\pi(C_\phi, \mc{L}|_{C_\phi})$ is of Arthur type, then it must belong to an Arthur packet $\Pi_\psi^\pure(G)$ for some Arthur parameter $\psi$ of $G$. 
By Vogan's conjecture $\Pi_{\phi_\psi}^\ABVpure(G)=\Pi_\psi^\pure(G)$.  
Therefore
$$\phi_\psi\in \Phi^{\ABV}(\pi(C_\phi, \mc{L}|_{C_\phi}))=\{\phi\},$$
and thus $\phi_\psi$ is equivalent to $\phi$, hence $\phi$ is Arthur type.  

Conversely, if $\phi$ is Arthur type, there exists an Arthur parameter $\psi$ of $G$ such that $\phi_\psi$ is equivalent to $\phi$.  
Assuming Vogan's conjecture,
$$\pi(C_\phi, \mc{L}|_{C_\phi})\in \Pi_{\phi_\psi}^\ABVpure(G)=\Pi_\psi^\pure(G),$$
and thus $\pi(C_\phi, \mc{L}|_{C_\phi})$ is Arthur type.

Finally, if $\psi'$ is some Arthur parameter for which $\pi(C_\phi, \mc{L}|_{C_\phi})\in \Pi_{\psi'}^\pure(G)$, then by Vogan's conjecture $\Pi_{\phi_{\psi'}}^\ABVpure(G)=\Pi_{\psi'}^\pure(G)$, and thus
$$\phi_{\psi'}\in \Phi^{\ABV}(\pi(C_\phi, \mc{L}|_{C_\phi}))=\{\phi_\psi\}.$$
Therefore $\psi$ and $\psi'$ are equivalent, hence $\Psi(\pi(C_\phi, \mc{L}|_{C_\phi}))=\{\psi\}$. 
\end{proof}

There are two special cases when $\bar C$ is smooth, namely, when $C$ is the open orbit or the closed orbit. 

In the closed orbit case, the result of Corollary \ref{cor: closed single} is closely related to the closure ordering conjecture, which says that for all $\pi\in \Pi_\psi(G)$, we have $C_{\phi_\psi}\subseteq \bar{C}_{\phi_\pi}$. See \cite{Xu}*{Conjecture 3.1}. The ABV-packet version of this conjecture is just \cite{CFMMX}*{Theorem 7.22 (b)}. The above conjecture was confirmed for $\Sp_{2n}$ and split $\SO_{2n+1}$ in \cite{HLLZ}. For the closed orbit $C=C_\phi$, if $\pi=\pi(C, \mc{L})\in \Pi_\psi(G)$, then by Xu's Conjecture, $C_{\phi_\psi}\subseteq \bar{C}_{\phi_\pi}=C$. Therefore, $\phi_\psi=\phi_\pi$ and $\Psi(\pi)=\{\psi\}$, which are consistent with Corollary \ref{cor: closed single} and Vogan's Conjecture \ref{con: vogans conjecture}.  

The case when $C$ is open will be discussed in more detail in the next section and we will see that Corollary \ref{cor: closed single} is exactly Enhanced Shahidi's conjecture, see Conjecture \ref{con: enhanced shahidi}. In this sense, Corollary \ref{cor: closed single} provides a generalization of Enhanced Shahidi's conjecture. For classical groups and their inner forms, Conjecture \ref{con: vogans conjecture} is settled in \cite{CHLLRXS}. Thus, Enhanced Shahidi's conjecture for these groups is now a theorem. We mention that for $\mathrm{Sp}_{2n}$ and split $\mathrm{SO}_{2n+1}$, Enhanced Shahidi's conjecture was proved in \cite{HLL} and \cite{HLLZ}. 

For more general groups other than classical groups, we expect that the results of Corollary \ref{cor: closed single} are unconditional. Namely, we expect that if $\phi$ is a Langlands parameter such that $\bar C_\phi$ is smooth, and $\mathcal{L}$ is an irreducible $H_\lambda$-equivariant local system on $\bar C_\phi$, then the representation $\pi(C_\phi,\mathcal{L}|_{C_\phi})$ is contained in at most one Arthur packet. In particular, if $\bar C_\phi$ is smooth, then $\pi(C_\phi, \1)$ is contained in at most one Arthur packet.

If $C$ is open, it is expected that $\pi(C,\1)$ is exactly the generic representation in $\Pi_{\lambda_\phi}$ for a fixed Whittaker datum, see \cite{CDFZ1}*{Conjecture 4.1}.
It is also interesting to ask if one can provide a representation-theoretic description of the class of representations of the form $\pi(C,\1)$, where $\bar{C}$ is smooth (or even more generally, representations of the form $\pi(C,\mathcal{L})$, where $\bar{C}$ is smooth and $\CL$ is an irreducible equivariant local system on $C$).  
At this moment, we are not sure if there are any non-open, non-closed orbits $C$ such that $\bar{C}$ is smooth and $\pi(C,\mc{L})$ is an Arthur type representation for classical groups.
Proving the non-existence or providing a classification of such representations (if there are any) would be interesting. 
We explore some examples in the next subsection. While for the exceptional group $G_2$, there are indeed Arthur type representations of the form $\pi(C,\1)$ for $\bar C$ smooth but neither open nor closed; see Example \ref{example-G2}.


\subsection{Smooth closure and Arthur-type representations}\label{sec:smooth closure}

In this section, we investigate various examples motivated by the following speculation:
    Given a classical group $G$, or $\GL_n$, if there is an Arthur parameter $\psi$ for which $\pi(C, \1_C)\in \Pi_\psi(G)$,  then $C$ is  open, closed, or $\bar{C}$ is singular. 

However, we will demonstrate below that it still holds for various infinitesimal parameters of $\GL_n$, despite not being classical.  

Given an Arthur parameter $\psi$, Vogan's conjecture \ref{con: vogans conjecture} predicts that
$$\pi(C_{\phi_\psi}, \1) \in \Pi_{\phi_{\psi}}^\ABV(G)=\Pi_\psi(G),$$
and thus, if our speculation holds, 
it must be that $C_\psi$ is open, closed, or has singular closure. 

In fact, we explore some examples of unramified infinitesimal parameters $\lambda$, for which, given a non-closed, non-open orbit $C\subseteq V_\lambda$, $\pi(C, \mathbbm{1})$ is of Arthur type if and only if $\bar{C}$ is not smooth.  
However, there are certainly cases for which $\bar{C}_\phi$ is singular, but $\phi$ is not of Arthur type.  
For example, in \cite{CFK}, an orbit $C_{\phi_{\text{KS}}}$ of a particular Vogan variety of $\GL_{16}$ is computed for which $\phi_{\text{KS}}$ is not of Arthur type. Meanwhile, the corresponding ABV-packet $\Pi_{\phi_{\text{KS}}}^\ABV(G)=\{\pi_{\text{KS}}, \pi_\psi\}$ for a particular Arthur type representation $\pi_\psi$. 

\begin{example}
Let $G$ be a split connected reductive group with a maximal torus of split rank $n$ and consider the infinitesimal parameter $\lambda$ of the Steinberg representation for $G$.  
Thus, for the simple roots $\Delta$,
$$V_\lambda\cong \bigoplus_{\alpha\in \Delta}\mf{g}_\alpha,$$
and $H_\lambda$ is the maximal torus corresponding to the choice of based root system.  
Every subset $S=\{\alpha_1, \ldots, \alpha_r\}$ of $\Delta$, determines an orbit
$$C_S:=\{a_1\alpha_1+\cdots +a_r\alpha_r \  | \ a_1, \ldots, a_r\in \bc^\ast\},$$
 and all orbits arise in this way. That is, the set of orbits is in bijection with the power set $\mc{P}(\Delta)$.  
 The closure of an orbit $C$ determined by any $r$ simple roots is isomorphic to $\mb{A}^r$, and hence smooth, and thus by Proposition \ref{prop: singleton}, for every Langlands parameter $\phi$,  $\Phi^\ABV(\pi(C_\phi,  \1))=\{\phi\}$.  

Each unramified Arthur parameter of $\GL_n$ is of the form 
$$\psi(w, x, y)=\bigoplus_{i=1}^m\chi_i\boxtimes S_{d_i}(x)\boxtimes S_{a_i}(y),$$
for some integers $a_i, d_i$ where $S_i$ denotes the $i$-dimensional irreducible representation of $\SL_2(\bc)$, and the $\chi_i$ denote unitary unramified characters of the Weil group $W_F$. 
Hence, by realizing each of the classical groups as subgroups of $\GL_n$, their Arthur parameters can also be realized in the above form.    

For $\GL_n$, the infinitesimal parameter of the Steinberg representation is equivalent to
$$\lambda(w)=\tx{diag}(\norm{w}^{(n-1)/2}, \ldots, \norm{w}^{(1-n)/2}).$$
The only Arthur parameters of the above form with infinitesimal parameter $\lambda$ are $\psi(w, x, y)=S_n(x)$ and $\psi(w, x, y)=S_n(y)$ corresponding to the open and closed orbits, respectively.  Since all irreducible equivariant local systems are trivial, the speculation holds for $\GL_n$.  As each Arthur parameter for the following groups considered must factor through $\GL_n$, an entirely similar argument tells us that the only orbits of Arthur type in the variety $V_\lambda$ where $\lambda$ is the infinitesimal parameter of the Steinberg representation, are the open and closed orbits. 

\end{example}

\begin{example}
For some natural number $n$ and each $i$ with $1\leq i \leq n$, consider the Arthur parameters
$$\psi_i(w, x, y):=S_2(x)^{\oplus i}\oplus S_2(y)^{\oplus(n-i)}.$$
The image of each $\psi_i$ can be taken to lie in $\GL_{2n}(\bc), \SO_{2n}(\bc)$ or $\Sp_{2n}(\bc)$ (for appropriate choice of $i$ in the latter cases), and thus is an Arthur parameter for $\GL_{2n}, \SO_{2n}$ and $\SO_{2n+1}$ respectively. 
In each case, the infinitesimal parameter is equivalent to 
$$\lambda(w)=\bigoplus_{i=1}^n(\norm{w}^{1/2}\oplus\norm{w}^{-1/2}),$$
where $\norm{\cdot}:W_F\to \bc^\times$ is the norm map, where $\norm{\mf{f}}=q$ and $\norm{\cdot}|_{I_F}=1$.  
Depending on the group, and the parity of $i$, $V_\lambda$ is isomorphic to a determinantal, symmetric determinantal, or skew-symmetric determinantal variety.  
In each case, each orbit is of the form $C_{\psi_i}$ for an appropriately chosen $i$. Hence, each orbit is of Arthur type, and the closure of each non-open, non-closed orbit has a singularity at the origin.  

The infinitesimal parameter of the Steinberg representation may be viewed as an extremal opposite of the infinitesimal parameter discussed above since the Steinberg infinitesimal parameter has all eigenvalues distinct, and the infinitesimal parameter above only has two distinct eigenvalues.  
(The Vogan variety of an unramified infinitesimal parameter with a single eigenvalue is a point.) 
The tables below highlight the relationship between orbits $C\subseteq V_\lambda$ of smooth closure and Arthur type representations for the Steinberg and 2-eigenvalue infinitesimal parameters for $\GL_{2n}$ and $\SO_n$.
We say that $\phi$ is of Arthur type if there exists an Arthur parameter such that $\phi$ is equivalent to $\phi_\psi$, and we say that a representation is of Arthur type if it belongs to an Arthur packet. 

\begin{center}
        \begin{tabular}{|c||c | c | c|}
        \hline 
        & & & \\
          Steinberg   & $\bar{C}_\phi$ Smooth &  $\phi$ Arthur type & $\pi(C_\phi, \1)$ Arthur Type \\
          \hline 
          \hline 
         $C_\phi$ Open/Closed  & Yes & Yes & Yes \\
         \hline 
         $C_\phi$ Non-Open/Closed & Yes & No & No \\
         \hline 
        \end{tabular}  
        \end{center}

        \begin{center}
        \begin{tabular}{|c||c | c |c |}
        \hline 
        & & & \\
          2 Eigenvalues   & $\bar{C}_\phi $ Smooth & $\phi$ Arthur type & $\pi(C, \1)$ Arthur Type \\
          \hline 
          \hline 
         Open/Closed  & Yes & Yes & Yes \\
         \hline 
         Non-Open/Closed & No & Yes  & Yes \\
         \hline 
        \end{tabular}
\end{center}
\end{example}

\begin{example}\label{ex:SO(7)}
In \cite{CFMMX}*{Chapter 16}, an infinitesimal parameter for $\SO_7$ which is neither of the above types is investigated.  
There are 8 corresponding Langlands parameters labeled as $\phi_i$ for $i\in \{0, 1, \ldots, 7
\}$, and other than the open and closed orbits, those of Arthur type are exactly those whose orbits do not have smooth closure. 

\begin{center}
\begin{tabular}{| c || c | c | c | c |}
\hline 
    & & & & \\
    L-parameter & $C_\phi$  Open/Closed & $\bar{C}_\phi$ Smooth  & $\phi$ Arthur & $\pi(C_\phi, \1)$ Arthur \\
    \hline 
    \hline 
    $\phi_0, \phi_7$  & Yes & Yes & Yes & Yes \\
    \hline 
    $\phi_2, \phi_4, \phi_5, \phi_6$ & No & No & Yes & Yes\\
    \hline 
    $\phi_1, \phi_3$ & No & Yes & No & No  \\
    \hline 
    \end{tabular}
\end{center}

Moreover, the corresponding Arthur packets are listed in \cite{CFMMX}*{Section 16.1}, and indeed, for $i\neq 0, 7$ (indices of the closed (resp. open) orbit) the representations $\pi(C_i, \mathbbm{1})$ of Arthur type are precisely those for which $\bar{C}$ is not smooth.  
Notice that $\bar{C}_3$ is smooth, and the representation $\pi(\phi_3, -)$ is of Arthur type, but is not the representation associated to the trivial sheaf.
\end{example}

It appears that this phenomena does not extend beyond $\GL_n$ and the classical groups, as seen in this example. 

\begin{example}\label{example-G2}
    In \cite{CFZ:unipotent}, every Vogan variety $V_\lambda$ was investigated for the group $G_2$ for unramified $\lambda$.
    For the infinitesimal parameter $\lambda_4$ appearing in \cite{CFZ:unipotent}*{Case 4, Section 1.3}, $V_{\lambda_4}$ has two Arthur type orbits with smooth closure, neither of which are open or closed. 
\end{example}


\section{The \texorpdfstring{$p$}{}-adic Kazhdan-Lusztig hypothesis and generic representations}\label{sec:pklh}

Throughout this section, we assume the Desiderata of Section \ref{llc} on the local Langlands correspondence and the $p$-adic Kazhdan-Lusztig hypothesis \ref{hypothesis}.  
The main result is Theorem \ref{thm: GP}, which (under the assumptions of this section) gives partial results on Conjecture \ref{conj-GP}. 
From the above, we prove Corollary \ref{cor: genABV} which is an analogue for ABV-packets. 
We also demonstrate in Corollary \ref{cor: VogSha} that Shahidi's enhanced genericity conjecture follows from Vogan's conjecture on ABV-packets (Conjecture \ref{con: vogans conjecture}).  

\subsection{Implications for standard representations}\label{sec:p-KLH}

\begin{proposition}\label{prop:maxform}
Let $G$ be a reductive $p$-adic group with infinitesimal parameter $\lambda$ satisfying the $p$-adic Kazhdan-Lusztig hypothesis. 
Then
\begin{enumerate}
\item For a pure inner form $[\delta]$ of $G$, let $C\subseteq V_\lambda$ be an $H_\lambda$-orbit such that for all $H_\lambda$-orbits $D\subseteq V_\lambda$ such that $C\subset\bar{D}, C\neq D$ and irreducible $\mc{F}\in \Loc_{H_\lambda}(D)$, $\pi(D, \mc{F})$ is not a representation of $G_\delta$.  Then, for every irreducible $\mc{L}\in \Loc_{H_\lambda}(C)$ for which $\pi(C, \mc{L})$ is a representation of $G_\delta$, the standard representation of $\pi(C, \mc{L})$ is irreducible and thus isomorphic to $\pi(C, \mc{L})$. 
\item Every standard representation for a representation in the pure L-packet of the open orbit is irreducible.
\item Every standard representation for a representation in the pure ABV-packet of the open orbit is irreducible.
\end{enumerate}
\end{proposition}
\begin{proof}
\begin{enumerate}
\item Let $C$ be as above and consider an irreducible $\mc{F}\in \Loc_{H_\lambda}(D)$ for some orbit $D$. 
By assumption, if $C\subseteq\bar{D}, C\neq D$ then $\pi(D, \mc{F})$ is not a representation of $G_\delta$, and therefore $[S(C, \mc{L}):\pi(D, \mc{F})]=0$.
  
By the Kazhdan-Lusztig hypothesis 
$$[S(C, \mc{L}):\pi(D, \mc{F})]=\sum_{n\in\bz}\dim\Hom_{\Loc_{H_\lambda}(C)}\lb\mc{H}^n(\mc{IC}(D, \mc{F})|_C, \mc{L}\rb,$$
but $\mc{IC}(D, \mc{F})$ is supported on $\bar{D}$.
Therefore, $[S(C, \mc{L}):\pi(D, \mc{F})]$ can only be non-zero if
$C\subseteq \bar{D}$, and hence is only nonzero for  $D=C$, in which case
\begin{align*}[S(C, \mc{L}):\pi(D, \mc{F})]&=\sum_{n\in\bz}\dim\Hom_{\Loc_{H_\lambda}(C)}\lb\mc{H}^n(\mc{IC}(C, \mc{F}))|_C, \mc{L}\rb\\
&=\sum_{n\in\bz}\dim\Hom_{\Loc_{H_\lambda}(C)}\lb\mc{H}^n(\mc{F}[\dim C]), \mc{L}\rb\\
&=\dim\Hom_{\Loc_{H_\lambda}(C)}\lb\mc{F}, \mc{L}\rb,
\end{align*}
which is 1 if $\mc{F}\cong \mc{L}$, and 0 otherwise.  
Therefore $S(C, \mc{L})\cong \pi(C, \mc{L})$ is irreducible. 
    \item As there are no orbits greater than the open orbit in the closure-ordering, this is a special case of (1).  
    \item  By \cite{CDFZ1}*{Theorem 5.1}, if $C_\phi$ is open,  
    $$\Pi_\phi^{\ABVpure}(G)=\Pi_\phi^{\pure}(G),$$
    and thus the result follows from (1). \qedhere
\end{enumerate}
\end{proof}

Recall that Lemma \ref{lem:irred_res} ensures that if $\mc{L}$ is an irreducible $H_\lambda$-equivariant local system on $\bar{C}$, then $\mc{L}|_C$ is irreducible. 

\begin{lemma}\label{lemma:smooth mult}
   Suppose that $C$ has smooth closure and $\mc{L}\in \Loc_{H_\lambda}(\bar{C})$ is irreducible.
   Then, for any orbit $D\subseteq \bar{C}$ and $\mc{F}\in \Loc_{H_\lambda}(D)$, 
   $$[S(D, \mc{F}):\pi(C, \mc{L}|_C)]=\dim \Hom(\mc{L}|_D, \mc{F}).$$
   Hence, if $D\subseteq \bar{C}$, and if $\mc{L}|_D$ is irreducible, there is a unique irreducible $\mc{F}\in \Loc_{H_\lambda}(D)$ for which $[S(D, \mc{F}):\pi(C, \mc{L}|_C)]\neq 0$, in which case $[S(D, \mc{F}):\pi(C, \mc{L}|_C)]=1$.
\end{lemma}
\begin{proof}
    As in the proof of \ref{lemma:memoirs}, for the inclusion $j:\bar{C}\hookrightarrow V$, we have $\mc{IC}(C, \mc{L}|_C)\cong j_!\mc{L}[\dim C]$.  
    Thus, by the Kazhdan-Lusztig hypothesis
    \begin{align*}
        [S(D, \mc{F}): \pi(C, \mc{L}|_C)]&=\sum_{n\in \bz}\dim \Hom_{\Loc_{H_\lambda}(D)}\lb\mc{H}^n(\mc{IC}(C, \mc{L}|_C)|_D, \mc{F}\rb\\
        &=\sum_{n\in \bz}\dim \Hom_{\Loc_{H_\lambda}(D)}\lb\mc{H}^n(j_!\mc{L})|_D, \mc{F}\rb\\
         &=\dim \Hom_{\Loc_{H_\lambda}(D)}\lb\mc{L}|_D, \mc{F}\rb.
    \end{align*}

    Since $\mc{L}|_D$ is irreducible, $[S(D, \mc{F}) : \pi(C, \mc{L}|_C)]$ is only non-zero for $\mc{F}\cong \mc{L}|_D$, in which case it is equal to 1. 
\end{proof}

\begin{proposition}\label{prop: generic appear in every standard module} 
Let $\lambda$ be an infinitesimal parameter of a reductive $p$-adic group for which the $p$-adic Kazhdan-Lusztig hypothesis holds.  
For an $H_\lambda$-orbit $C\subseteq V_\lambda$,
\begin{enumerate}
\item[ $(1)$] If $C$ is open, then for every irreducible $\mc{L}\in \Loc_{H_\lambda}(C)$, $\pi(C, \mc{L})$ is isomorphic to its own standard representation.  
 \item[ $(2)$] If $C$ has smooth closure, then for each $H_\lambda$-orbit $D\leq C$, $\pi(C, \1)$ appears in the standard representation of $\pi(D, \1)$ with multiplicity 1, and in no other standard representations. 
 \item[ $(3)$] $C$ is the open orbit if and only if $\pi(C, \1)$ is isomorphic to its own standard representation.  
\end{enumerate}
\end{proposition}

\begin{proof}
\begin{enumerate}
\item As there are no orbits greater than the open orbit in the closure ordering, the result follows by Proposition \ref{prop:maxform} (1). 
\item By Lemma \ref{lemma:smooth mult},
\begin{align*}
    [S(D, \1): \pi(C, \1)]&=\dim \Hom(\1_{\bar{C}}|_D, \1_D)\\
    &=\begin{cases}
    \dim\Hom(\1_D, \1_D), & D\leq C\\
    0, & \tx{otherwise}
    \end{cases}\\
    &=\begin{cases}
    1, & D\leq C\\
    0, &\tx{otherwise.}
    \end{cases}
\end{align*}
\item If $C$ is the open orbit, then this follows from (1). 

Conversely, suppose that $\pi(C, \1)$ is isomorphic to its own standard. 
If $D$ is the open orbit, then $\bar{D}=V_\lambda$ is smooth, hence by (2) 
$$[\pi(C, \1):\pi(D, \1)]=[S(C, \1):\pi(D, \1)]=1,$$
but this is only possible if $\pi(C, \1)\cong \pi(D, \1)$, which is to say $C=D$ is open. \qedhere
\end{enumerate}
\end{proof}

Following the main result of the next section, we will give an interpretation of Proposition \ref{prop: generic appear in every standard module} (2) in terms of generic representations. 

\begin{remark}
Heiermann \cite{H} proved that if a Vogan $L$-packet $\Pi_\phi^{\mathrm{pure}}(G)$ of a quasi-split classical group contains a generic representation, then for each $\pi\in \Pi_{\phi}^{\mathrm{pure}}(G)$, the standard module of $\pi$ is irreducible and hence isomorphic to $\pi$. Conjecture \ref{conj-GP} (which is known for classical groups) says that $\Pi_\phi^{\mathrm{pure}}(G)$ is generic if and only if $C_\phi$ is open, see \cite{CDFZ1}* {Conjecture 4.1} and Theorem \ref{thm: GP}. 
Thus, Proposition \ref{prop: generic appear in every standard module} offers a geometric way to state Heiermann's result and extend it to non-classical groups satisfying the $p$-adic Kazhdan-Lusztig hypothesis. 
\end{remark}

\subsection{Generic \texorpdfstring{$L$}{}-packets}\label{sec:genL}

Recall that Conjecture \ref{conj: genericL intro}, which is equivalent to the Gross-Prasad Conjecture \ref{conj-GP}, predicts that $\Pi_\phi(G)$ contains a generic representation if and only if $C_\phi$ is open. 
The result of \cite{CDFZ1}*{Theorem 4.1} is that the above holds for quasi-split classical groups.  
In this section, we will prove Theorem \ref{thm: GP}:  If $\phi$ is a Langlands parameter such that the $p$-adic Kazhdan-Lusztig hypothesis is known for $\lambda_\phi$, then the forward direction of Conjecture \ref{conj: genericL intro} holds. 
Moreover, if there exists a generic representation with infinitesimal parameter $\lambda$, then the converse holds as well.

For example, if $G$ is unramified, then by \cite{CasSha}, for every unramified character $\chi$ of a maximal split torus in $G$, the representation $I_B^G(\chi)$ contains a generic representation.  
Choosing an element $w$ of the Weyl group of $G(F)$ such that $\chi^w$ is $B$-positive, $I_B^G(\chi^w)$ is a standard representation.  Thus, by the $p$-adic Kazhdan-Lusztig hypothesis, every irreducible subquotient must have the same infinitesimal parameter.  
Therefore, if $\lambda$ is the infinitesimal of an unramified group corresponding to a subquotient of $I_B^G(\chi)$, then there is a generic representation having infinitesimal parameter $\lambda$.

The cases of Conjecture \ref{conj: genericL intro} which are covered in \cite{CDFZ1} rely on \cite{GI}, which imposes some (mildly) stronger conditions on the local Langlands correspondence than the proof we give here.
On the one hand, the proof given here does not depend on \cite{GI}, and thus requires weaker assumptions on the local Langlands correspondence. On the other hand, we impose the stronger assumption of the $p$-adic Kazhdan-Lusztig hypothesis, although it is already known in many cases, as discussed following the statement of Hypothesis \ref{hypothesis}.  
As a consequence, in the next section, we will be able to strengthen \cite{CDFZ1}*{Proposition 4.4}, and using this result, we demonstrate that if $G$ is a group for which Arthur's conjectures (as in \cite{Arthur:book}) are known, as well as Vogan's Conjecture \ref{con: vogans conjecture} holds, then Shahidi's enhanced genericity conjecture also holds.

\begin{theorem}\label{thm: GP}
Let $\lambda$ be an infinitesimal parameter of a quasi-split reductive group $G$ for which the $p$-adic Kazhdan-Lusztig hypothesis holds.   
If $\phi$ is a Langlands parameter with infinitesimal parameter $\lambda$ such that $\Pi_\phi(G)$ contains a generic representation, then $C_\phi$ is open. 

If \ $\Pi_\lambda(G)$ contains a generic representation, then if $C_\phi$ is open, $\Pi_\phi(G)$ contains a generic representation. 
\end{theorem}
\begin{proof}
Suppose $\pi\in \Pi_\phi(G)$ is $\mf{w}$-generic. Labeling the representations according to the corresponding Whittaker normalization $J(\mf{w})$ as in Definition \ref{piCL}, $\pi\cong \pi(C, \1)$.  By the standard representation conjecture \cite{CS}, which was proved in \cite{opdamh}, $\pi$ is isomorphic to its own standard representation, and thus by Proposition \ref{prop: generic appear in every standard module} (3) $C_\phi$ is open.

Suppose $\pi\in \Pi_\lambda(G)$ is generic, and $C_\phi$ is open.     
Letting $\mf{w}$ be the Whittaker data for $\pi$, and choosing the normalization with respect to $\mf{w}$, $\pi\cong \pi(C_{\phi_\pi}, \1)$.  
 By the forward direction of this proof $C_{\phi_\pi}$ must be the open orbit $C_\phi$, and thus $\Pi_\phi(G)=\Pi_{\phi_\pi}(G)$ contains a generic representation. 
 
\end{proof}

Thus, we can see that the converse direction of the above statement follows, in general, from the forward direction, and the claim that for every infinitesimal parameter $\lambda:W_F\to {}^LG$ there exists an irreducible generic representation $\pi$ such that $\lambda_\pi$ is equivalent to $\lambda$. 

We can now give a more representation-theoretic description of Proposition \ref{prop: generic appear in every standard module} (2):

\begin{corollary}
    Let $\pi$ be a $\mf{w}$-generic representation of a quasi-split group $G$ and let $\phi$ be its Langlands parameter.  
    Then, with respect to the Whittaker normalization $J(\mf{w})$, for every orbit $D$ of $V_{\lambda_\phi}$, and an irreducible $\mc{L}\in \Loc_{H_\lambda}(D)$,
    $$[S_{\pi(D, \mc{L})}: \pi] =\begin{cases}
        1, & \mc{L} =\1 \\
        0, & \tx{else}
    \end{cases}$$
\end{corollary}
\begin{proof}
    By our choice of Whittaker normalization, $\pi\cong \pi(C_\phi, \1)$. 
    Hence, by Theorem \ref{thm: GP}, we know that $C_\phi$ is open. 
    As every orbit $D$ satisfies $D\leq C_\phi$, the result follows from
    Proposition \ref{prop: generic appear in every standard module} (2).
\end{proof}

We also demonstrate an easy Corollary to Theorem \ref{thm: GP} for generic supercuspidal representations.  

\begin{corollary}
Let $\lambda$ be an infinitesimal parameter for a quasi-split reductive $p$-adic group $G$ for which the $p$-adic Kazhdan-Lusztig hypothesis \ref{hypothesis} and the desiderata of the local Langlands correspondence in Section \ref{llc} hold. Let $\pi$ be a supercuspidal generic representation of $G(F)$, and $\phi=\phi_\pi: W_F\times \SL_2(\C)\to \Lgroup{G}$ be its local Langlands parameter and $\lambda=\lambda_\phi$. Then $V_\lambda=\{0\}$ and $\phi_\pi|_{\SL_2(\C)}=1$.
\end{corollary}

\begin{proof}
The proof is along the same lines as the proof of Theorem~\ref{thm: GP}. Note that $V_\lambda=\{0\}$ implies $ \phi_\pi|_{\SL_2(\C)}=1$, thus it suffices to show that $V_\lambda=\{0\}.$ Assume that $C_0$ and $C$ are the respective closed and open orbits of $V_\lambda$.  Since $\pi$ is assumed to be generic, then $x_\phi\in C$ by Theorem~\ref{thm: GP}. 
Desiderata (2) in  \cite{CDFZ1}*{Section 1.7} implies that we can choose a Whittaker normalization such that $\pi\cong \pi(C,\1)$. Now Proposition \ref{prop: generic appear in every standard module} implies that $[S_{\pi(C_0, \1)} :\pi]=1$.
Therefore, $\pi\cong \pi(C, \1)$ and $\pi(C_0, \1)$ have the same cuspidal support, but since $\pi$ is assumed to be supercuspidal, this is only possible if $C=C_0$ and thus $V_\lambda= \{ 0\}.$
\end{proof}


\subsection{Generic A- and ABV-packets}\label{sec:genABV}

In this section, we will use Theorem
\ref{thm: GP} to strengthen the result
\cite{CDFZ1}*{Proposition 4.4} and make further progress to Shahidi's enhanced genericity conjecture.  

Recall that a Langlands parameter $\phi$ is said to be tempered if $\phi|_{W_F}$ has bounded image, and we say that an Arthur parameter $\psi$ is tempered if $\phi_\psi$ is tempered. 
The well-known conjecture \cite{Shahidi: plancherel}*{Conjecture 9.4} predicts that if $\phi$ is tempered, then $\Pi_\phi(G)$ contains a generic representation. 
More recently, Shahidi proposed an enhanced version of this conjecture.  

\begin{conjecture}[\cite{LLS}*{Conjecture 1.2}]\label{con: enhanced shahidi}
An Arthur packet $\Pi_\psi(G)$ contains a generic representation if and only if $\psi$ is tempered.  
\end{conjecture}

Inspired by Shahidi's enhanced genericity conjecture, an analogue for ABV packets was proved in \cite{CDFZ1}*{Theorem 4.4} under the assumption of Theorem \ref{thm: GP}.  Thus, assuming the $p$-adic Kazhdan-Lusztig hypothesis, we generalize one direction to arbitrary quasi-split connected reductive groups.

   \begin{corollary}\label{cor: genABV}
    Let $G$ be a quasi-split reductive group over a $p$-adic field $F$, and let 
    $\lambda$ be an infinitesimal parameter for which the $p$-adic Kazhdan-Lusztig hypothesis holds. 
    If \ $\Pi_\phi^\ABVpure(G)$ contains a generic representation, then $C_\phi$ is open.  

    If \ $\Pi_\lambda(G)$ contains a generic representation and $C_\phi$  is open, then $\Pi_\phi^{\ABVpure}(G)$ contains a generic representation. 
\end{corollary}

\begin{proof}
   Let $\pi$ be a generic representation in $\Pi_\phi^\ABVpure(G)$.  
    Letting $\lambda=\lambda_\phi$, and $\phi_\pi$ the Langlands parameter of $\pi$,
    by Theorem \ref{thm: GP} $C_{\phi_\pi}$ is open, and $\pi\cong \pi(C_\pi, \1)$.  
    As $\Pi_{\phi_\pi}(G)\subseteq \Pi_{\phi_\pi}^\ABVpure(G)$, we know that $\pi(C_{\phi_\pi}, \1)\in \Pi_{\phi_\pi}^\ABVpure(G)$, but by Proposition \ref{prop: singleton} $\pi(C_{\pi}, \1)$ belongs to a single ABV-packet. Therefore, $\phi$ and $\phi_\pi$ are equivalent. 

 If $\Pi_\lambda(G)$ contains a generic representation, then by the converse direction of Theorem \ref{thm: GP} if $C_\phi$ is open, then $\Pi_\phi(G)$ contains a generic representation, and since
 $$\Pi_\phi(G)\subseteq \Pi_\phi^\ABVpure(G),$$
  the result follows. 
\end{proof}

The result \cite{CDFZ1}*{Proposition 4.8} states that for a classical quasi-split group $G$, if one containment $\Pi_\psi^\pure(G)\subseteq \Pi_{\phi_\psi}^\ABVpure(G)$ of Vogan's Conjecture \ref{con: vogans conjecture} holds, then Shahidi's enhanced genericity conjecture also holds.  
One part of the proof in \cite{CDFZ1} that depends on $G$ being classical is \cite{CDFZ1}*{Proposition 4.7}, which we now generalize.  

\begin{corollary}\label{corollary: classical}
Let $G$ be a quasi-split reductive group over $F$, and assume the Kazhdan-Lusztig hypothesis \ref{hypothesis}. If $L(s,\phi,\Ad)$ has a pole at $s=1$, then $\Pi_\phi^{\ABVpure}(G)$ does not contain a generic representation.
\end{corollary}
\begin{proof}
   Suppose $\pi\in \Pi_\phi^\ABVpure(G)$ is generic.  
    By Corollary \ref{cor: genABV} $C_\phi$ must be open. Hence, by \cite[Proposition 3.5]{CDFZ1}  $L(s, \phi, \Ad)$ is regular at $s=1$.
\end{proof}

The only other part of the proof of \cite{CDFZ1}*{Proposition 4.8} that technically assumes that $G$ is classical is the fact that it relies on Arthur's conjectures on A-packets.  
Thus, we can upgrade the result to a statement conditional on Arthur's conjectures.  

\begin{corollary}\label{cor: VogSha}
    Let $G$ be a quasi-split classical group for which Arthur's conjectures, as stated in \cite{Arthur:book}, on local $A$-packets are known, and  assume Vogan's Conjecture \ref{con: vogans conjecture}: $\Pi_\psi^\pure(G)= \Pi_{\phi_\psi}^\ABVpure(G)$. 
    Then,
    \begin{enumerate}
        \item The forward direction of Shahidi's enhanced Conjecture \ref{con: enhanced shahidi} holds: If $\Pi_\psi(G)$ contains a generic representation, then $\phi_\psi$ is tempered. 
        \item If $\phi_\psi$ is tempered, and $\Pi_{\lambda_\psi}(G)$ contains a generic representation, then $\Pi_\psi(G)$ contains a generic representation. 
    \end{enumerate}
\end{corollary}
\begin{proof}
The proof that follows is similar to the one given in \cite{CDFZ1}*{Proposition 4.8}.  
\begin{enumerate} 
\item Suppose that $\pi\in \Pi_\psi^\pure(G)$ is generic.  
By Vogan's conjecture
$$\Pi_\psi^\pure(G)= \Pi_{\phi_\psi}^\ABVpure(G),$$
hence $\Pi_{\phi_\psi}^\ABVpure(G)$ contains a generic representation. Hence, by Corollary \ref{cor: genABV}, $C_{\phi_\psi}$ is open. 
By \cite{CDFZ1}*{Proposition 3.2}, it must be that $\phi_\psi$ is tempered. 

 \item If $\phi_\psi$ is tempered, then by \cite[Proposition 3.2]{CDFZ1} $C_{\phi_\psi}$ is open.  
 Since $\Pi_{\lambda_\psi}(G)$ contains a generic representation, by Corollary \ref{cor: genABV} $\Pi_{\phi_\psi}^{\ABVpure}(G)$ contains a generic representation, and thus by Vogan's conjecture, so does $\Pi_\psi^\pure(G)$.   
 \end{enumerate}
\end{proof}



\end{document}